\documentclass[12pt,twoside]{extarticle}
\usepackage{mathtext}
\usepackage{titlesec}
\usepackage[cp1251]{inputenc}
\usepackage[T2A]{fontenc}
\usepackage[russian,english]{babel}
\usepackage{fancyhdr}
\usepackage{amsthm}
\usepackage{dsfont}
\usepackage{indentfirst}
\usepackage{mathrsfs}
\usepackage{array}
\usepackage{amssymb}
\usepackage{amsmath}
\usepackage{eucal}
\usepackage{dubloper}
\usepackage{bm}
\usepackage{epsfig}
\usepackage{wrapfig}

\titleformat{\section}[block]{\large\bfseries\filcenter}{\large\bfseries\thesection. }{0pt}{}
\titleformat{\subsection}[block]{\bfseries\filcenter}{\bfseries\thesubsection. }{0pt}{}
\titleformat{\subsubsection}[block]{\bfseries\filcenter}{\bfseries\thesubsubsection. }{0pt}{}

\newcounter{mytheorem}[section]

\edef\lim{\lim\limits}
\renewcommand{\frac}[2]{\dfrac{#1}{#2}}
\numberwithin{equation}{section}
\numberwithin{figure}{section}
\begin{document}

\setlength{\abovedisplayskip}{2mm} 
\setlength{\belowdisplayskip}{2mm} 
\parindent=8,8mm
\renewcommand{\headrulewidth}{0pt}
\edef\sum{\sum\limits}

\thispagestyle{empty}

\begin{center}

\

\

\

\

\textbf{\Large On finite groups whose}

\textbf{\Large Sylow subgroups are submodular}

\bigskip

\textbf{\bf V.\,A.\,Vasilyev}

\end{center}

\bigskip

\begin{center}
{\bf Abstract}
\end{center}

\medskip

\noindent\begin{tabular}{m{8mm}m{138mm}}

&{\small~~~
A subgroup $H$ of a finite group $G$ is called submodular in $G$,
if we can connect $H$ with $G$ by a chain of
subgroups, each of which is modular (in the sense of Kurosh) in the next.
If a group $G$ is supersoluble and every Sylow subgroup of $G$ is submodular in $G$, then
$G$ is called strongly supersoluble. The properties of groups with submodular
Sylow subgroups are obtained. In particular, we proved that in a group every
Sylow subgroup is submodular if and only if the group is Ore dispersive
and every its biprimary subgroup is strongly supersoluble.

{\bf Keywords:} finite group, modular subgroup, submodular subgroup,
strongly supersoluble group, Ore dispersive group.}

\smallskip

Mathematics Subject Classification (2010): 20D10, 20D20, 20D40. 


\end{tabular}

\bigskip

{\parindent=0mm
\textbf{\bf {\large Introduction}}
}
\bigskip

Throughout this paper, all groups are finite. The notion of a normal subgroup takes a central place in the theory of groups. One of its generalizations is the notion of a modular subgroup, i.e. a modular element (in the sense of Kurosh \cite[Chapter~2, p.~43]{sch2}) of a lattice of all subgroups of a group.
Recall that a subgroup $M$ of a group $G$ is called modular in $G$,
if the following hold:

1) $\langle X,M \cap Z \rangle=\langle X,M \rangle \cap Z$ for all $X \leq G, Z \leq G$ such that
$X \leq Z$, and

2) $\langle M,Y \cap Z \rangle=\langle M,Y \rangle \cap Z$ for all $Y \leq 
G, Z \leq G$ such that
$M \leq Z$.

Properties of modular subgroups were studied in the book \cite{sch2}.
Groups with all subgroups are modular were studied by R.~Schmidt \cite{sch2},
\cite{sch11} and
I.~Zimmermann \cite{zim}. By parity of reasoning with subnormal subgroup, in
\cite{zim} the notion of a submodular subgroup was introduced.

{\bf Definition \cite{zim}.} {\sl A subgroup $H$ of a group $G$ is called
submodular in $G$, if
there exists a chain of subgroups
$H=H_0 \leq H_1 \leq \ldots \leq \ H_{s-1} \leq H_{s} = G$ such that $H_{i-1}$ is a
modular subgroup in $H_{i}$ for $i=1,\ldots,s$.}

If $H \neq G$, then the chain can be compacted to maximal modular subgroups.

It's well known that in a nilpotent group every Sylow subgroup is normal (subnormal).
In the paper \cite{zim} groups
with submodular subgroups were studied. In particular, 
it was proved that in a supersoluble group $G$ every Sylow subgroup is submodular
if and only if $G/F(G)$ is abelian
of squarefree exponent.
A criterion of the submodularity of Sylow subgroups in an arbitrary group was found.

This paper is devoted to the further study of groups with submodular Sylow subgroups.

A group we will call strongly supersoluble, if it is supersoluble and every its
Sylow subgroup is submodular in it.

The class of all strongly supersoluble groups we will denote $s\mathfrak{U}$.
We proved that the class of groups $s\mathfrak{U}$ is a hereditary
local formation. We obtained that a group is strongly supersoluble if and only if it is
metanilpotent and every its Sylow subgroup is submodular. The class of all groups
with submodular Sylow subgroups we denote $sm\mathfrak{U}$. We proved that
$sm\mathfrak{U}$
forms a hereditary local formation and its local screen was found.
We established that in a group every Sylow subgroup is submodular if and only if
the group is Ore dispersive and every its biprimary subgroup is strongly supersoluble.



\bigskip

{\parindent=0mm
\textbf{\bf {\large 1. Preliminaries}}
}

\bigskip

We use the standard notation and terminology (see \cite{shem} and \cite{DH}).
Recall some of them.

Let $G$ be a group.
${\rm {Syl}}_p(G)$ is a set of all Sylow $p$-subgroups of $G$
for some prime $p$;
${\rm {Syl}}(G)$ is a set of all Sylow subgroups of $G$;
$M_G$ is the core of subgroup
$M$ of $G$, i.e. the intersection of all subgroups conjugated with $M$ in $G$;
$F(G)$ is the Fitting subgroup of $G$, i.e. the product of all normal
nilpotent subgroups of $G$;
$F_p(G)$ is a $p$-nilpotent radical of $G$, i.e. the product of all normal $p$-nilpotent subgroups of $G$,
$p$ is some prime.

A group $G$ of order $p_1^{n_1}p_2^{n_2}\dots p_n^{n_k}$ is called Ore dispersive \cite[p.~251]{shem},
if $p_1>p_2>\dots >p_n$ and
$G$ has a normal subgroup of order $p_1^{n_1}p_2^{n_2}\dots p_n^{n_i}$ for every $i=1, 2,\dots, k$.

We use the following notation  for concrete classes of group:
$\frak{S}$ is the class of all soluble groups;
$\frak{U}$ is the class of all supersoluble groups;
$\frak{N}$ is the class of all nilpotent groups;
$\frak{A}(p-1)$ is the class of all abelian groups of exponent dividing $p-1$.
By ${\cal M}(\frak{X})$ is denoted the class of all minimal non $\frak{X}$-groups,
i.e. such groups $G$ for which all proper subgroups of $G$ are contained in $\frak{X}$, but $G\not\in \frak{X}$.

A class of groups $\frak F$ is called a {\it formation} if the following conditions hold:
(a) every quotient group of a group lying in $\frak F$ also lies in
$\frak F$; (b) if $H/A\in \frak F$ and $H/B\in \frak F$ then
$H/A\cap B\in \frak F$.

A formation $\frak{F}$ is called
{\it hereditary} whenever $\frak{F}$ together with every group contains
all its subgroups, and {\it saturated}, if $G/\Phi(G)\in\frak{F}$ implies that $G\in\frak{F}$.

The $\frak F$-residual of a group $G$ for nonempty formation is denoted by $G^{\frak {F}}$,
i.e. the smallest normal subgroup of $G$ with $G/G^{\frak {F}}\in
\frak {F}$.

A function $f:\Bbb{P}\rightarrow \{$formations$\}$ is called a {\it local screen}.
A formation $\frak{F}$ is called {\it local},
if there exists a local screen $f$ such that $\frak{F}$
coincides with the class of groups $(G | G/C_G(H/K)\in f(p)$
for every chief factor $H/K$ of $G$ and
$p\in\pi(H/K))$.
It's denoted by $\frak{F}=LF(f)$.

Recall that a subgroup $H$ of $G$ is called maximal modular in $G$,
if $H$ is modular in $G$ and from $H \leq M < G$ it always follows $H = M$ for every modular subgroup $M$ in $G$.

{\bf Lemma 1.1 \cite[Lemma 1]{zim}.} {\sl Let $G$ be a group and $T \leq G$.
Then the following hold:

$1)$ if $T$ is submodular in $G$ and $U \leq G$, then $U \cap T$
is submodular in $U$;

$2)$ if $T$ is submodular in $G$, $N$ is normal in $G$ and $N \leq T$, then
$T/N$ is submodular in $G/N$;

$3)$ if $T/N$ is submodular in $G/N$, then $T$ is submodular in $G$;

$4)$ if $T$ is submodular in $G$, then $T^x$ is submodular in $G$ for every
$x \in G$;

$5)$ if $T_1$ and $T_2$ are submodular subgroups in $G$, then $T_1 \cap T_2$
is a submodular subgroup in $G$;

$6)$ if $T$ is submodular in $G$, then $TN$ is submodular in $G$
for every normal in $G$ subgroup $N$.

}

\medskip

{\bf Lemma 1.2 \cite[Lemma 1]{sch11}.} {\sl A subgroup $M$ of a group $G$ is maximal modular in $G$ if and only if either $M$ is maximal normal subgroup in $G$, or $G/M_G$ is nonabelian of order $pq$, where $p$ and $q$ are primes.
}

\medskip

{\bf Lemma 1.3 \cite[Lemma 2]{hein}.} {\sl Let $G=AB$ be a product of nilpotent subgroups $A$ and $B$, and $G$ has a minimal normal subgroup $N$ such that $N=C_G(N) \neq G$. Then

$1)$ $A \cap B = 1$;

$2)$ $N \subseteq A \cup B$;

$3)$ if $N \leq A$, then $A$ is a $p$-group for some prime $p$ and $B$ is $p'$-group.

}

\medskip

{\bf Lemma 1.4 \cite[Lemma 3.9 (1)]{shem}.} {\sl If $H/K$ is a chief factor of a group $G$ and $p \in \pi(H/K)$, then $G/C_G(H/K)$ is not contain nonidentity
normal $p$-subgroups, and besides $F_p(G) \leq C_G(H/K)$.}

\medskip

{\bf Theorem 1.5 \cite[Theorem 1.4]{wei}.} {\sl Let $H/K$ be a $p$-chief
factor of a group $G$. Then $|H/K|=p$ if and only if $Aut_G(H/K)$ is abelian of exponent dividing $p-1$.}

\medskip

{\bf Theorem 1.6 \cite[Chapter IV, Theorem 4.6]{DH}.} {\sl A formation is saturated if and only if it is local.}

{\bf Lemma 1.7 \cite[Lemma 4.5]{shem}.} {\sl Let $f$ be a local screen of a formation $\frak{F}$. A group $G$ belongs to $\frak{F}$ if and only if $G/F_p(G)\in f(p)$ for every $p\in\pi(G)$.}

We need the following property of the class of all supersoluble groups (see,
for example, \cite[p.~35]{shem} or \cite[p.~358]{DH}).

{\bf Lemma 1.8.} {\sl The class of all supersoluble groups has a local screen $f$ such that $f(p)=\frak{A}(p-1)$ for every prime $p$.}

\bigskip

{\parindent=0mm
\textbf{\bf {\large 2. Strongly supersoluble groups}}
}

\bigskip

{\bf Lemma 2.1.} {\sl Let $p$ be the largest prime divisor of $|G|$ and
$G_p \in \text{\rm Syl}_p(G)$. If $G_p$ is submodular subgroup in $G$,
then $G_p$ is normal in $G$.}

{\bf Proof.} We will use an induction by $|G|$.
We can consider that $G_p \neq G$ and there exists a chain of subgroups
$G_p=H_0 < H_1 < \ldots < H_{s-1} < H_{s} = G$ such that $H_{i-1}$ is maximal modular subgroup in $H_{i}$ for $i=1,\ldots,s$. By induction, $G_p$ is
normal in $H_{s-1} = M$. By Lemma 1.2, either $M$ is normal in $G$, or
$G/M_G$ is nonabelian of order $rq$, where $r$ and $q$ are primes. In the first case,
$G_p$ is normal in $G$. So let $|G/M_G|=rq$ and $G/M_G$ be a nonabelian group. It follows $|G:M|$ is a prime different from $p$. So we can assume that $|G:M|=q \neq p$.
If $N_G(G_p) \neq G$, then by the Theorem of Sylow $|G:M| = |G:N_G(G_p)|
= q \equiv 1\ ({\rm{mod}}\ p)$. We got a contradiction with $q < p$.
So $N_G(G_p) = G$. Lemma is proved.

\medskip

{\bf Corollary 2.1.1 \cite[Proposition 9]{zim}.} {\sl If every Sylow subgroup 
of the group $G$ is submodular in $G$, then $G$ is Ore dispersive.}

\medskip

{\bf Definition 2.2.} {\sl A group $G$ we will call strongly supersoluble if $G$ is supersoluble and every Sylow subgroup of $G$ is submodular in $G$.}

Denote $s\mathfrak{U}$ the class of all strongly supersoluble groups.

\medskip

{\bf Proposition 2.3 \cite[Proposition 10]{zim}.} {\sl A group $G$ is strongly
supersoluble if and only if $G$ is supersoluble and $G/F(G)$ is abelian
of squarefree exponent.}

\medskip

 In the paper we denote $\mathfrak{B}$ the class of all abelian groups of exponent free from squares of primes.

{\bf Lemma 2.4.} {\sl The class of groups $\mathfrak{B}$
is a hereditary formation.}

{\bf Proof.} It's clear, if $G \in \mathfrak{B}$, then
$H \in \mathfrak{B}$ and $G/N \in \mathfrak{B}$ for any subgroup $H$
and any normal subgroup $N$ of $G$.

Let $G$ be a group of the smallest order such that $G/N_i \in \mathfrak{B},
N_i \unlhd G, i=1,2$, but $G/N_1 \cap N_2 \not \in \mathfrak{B}$.

 If $N_1 \cap N_2 \neq 1$, then in $N_1 \cap N_2$ there is a normal subgroup $K$ in $G$ .
 From $|G/K| < |G|$ and $G/K/N_i/K \simeq G/N_i \in \mathfrak{B}$, $i=1,2$, it follows $G/K/(N_1/K \cap N_2/K) \simeq G/N_1 \cap N_2 \in \mathfrak{B}$. This contradicts the choice of $G$.

 Let $N_1 \cap N_2 = 1$. Since $ \mathfrak{A}$ is a formation and $G/N_i \in \mathfrak{B}
 \subseteq \mathfrak{A}$, then $G \in  \mathfrak{A}$. Let's show that the exponent of $G$ is free from squares of primes. Let $z$ be an element of order $q^n$ from $G$, where
$q$ is a prime, and $Z = \langle z \rangle$.
Assume that $Z \cap N_1 \neq 1$ and $Z \cap N_2 \neq 1$. Since $Z$ is cyclic $q$-group, then there exists $i \in \{1,2\}$ such that $Z \cap N_i \leq Z \cap N_{3-i}$.
 Hence we get a contradiction $1 \ne Z \cap N_i \leq (Z \cap N_1) \cap (Z \cap N_2) =1$.
Hence $Z \cap N_j = 1$ for some $j \in \{1,2\}$. Then from $G/N_j \in
\mathfrak{B}$ and $ZN_j/N_j \simeq Z$ it follows $n < 2$ and $G \in \mathfrak{B}$.
Lemma is proved.

\medskip

Note that the class of groups $\mathfrak{B}$ is not saturated. For example,
a cyclic group $G=\langle z | z^4=1\rangle \not \in \mathfrak{B}$,
but $G/\Phi(G) \in \mathfrak{B}$.

\medskip

{\bf Theorem 2.5.} {\sl Let $G$ be a group. Then the following hold:

$1)$ if $G$ is strongly supersoluble, then every subgroup of $G$ is strongly supersoluble;

$2)$ if $G$ is strongly supersoluble and $N \unlhd G$, then $G/N$ is strongly supersoluble;

$3)$ if $N_i \unlhd G$ and $G/N_i$ is strongly supersoluble for $i=1,2$,
then $G/N_1 \cap N_2$ is strongly supersoluble;

$4)$ if $H_i \unlhd G$, 
$H_i$ is strongly supersoluble, $i=1,2$ and $H_1 \cap H_2 =1$, then
$H_1 \times H_2$ is strongly supersoluble;

$5)$ if $G/\Phi(G)$ is strongly supersoluble, then $G$ is strongly supersoluble;

$6)$ the class of groups $s\mathfrak{U}$ is a hereditary saturated formation.
}

{\bf Proof.} Show the validity of 1). Let $G\in s\mathfrak{U}$ and $H\leq G$. In view of Propositon 2.3 and Lemma 2.4, from $G/F(G)\in \mathfrak{B}$ it follows that $H/H\cap F(G)\simeq HF(G)/F(G)
\in \mathfrak{B}$. Since $\mathfrak{B}$ is a homomorph and $H\cap F(G)\leq F(H)$,
we get $H/F(H)\simeq H/H\cap F(G)/F(H)/H\cap F(G)\in \mathfrak{B}$. By Proposition 2.3,
$H\in s\mathfrak{U}$.

Prove Statement 2). Let $G\in s\mathfrak{U}$ and $N \unlhd G$.
By Proposition 2.3,
$G/F(G)\in \mathfrak{B}$. Since $\mathfrak{B}$
is a homomorph, from $F(G)N/N\leq F(G/N)=F/N$ we conclude
$G/N/F(G/N)\simeq G/F\simeq G/F(G)/F/F(G)\in \mathfrak{B}$. By Proposition 2.3, $G/N\in s\mathfrak{U}$.

Prove Statement 3). Let $G$ be a group of the smallest order such that
$N_i \unlhd G$ and $G/N_i\in s\mathfrak{U}$ for $i=1,2$, but
$G/N_1 \cap N_2\not\in s\mathfrak{U}$. Since $\mathfrak{U}$ is a formation, $G/N_1 \cap N_2\in \mathfrak{U}$.

If $N_1\cap N_2\not=1$, then take from $N_1\cap N_2$ a subgroup $K\unlhd G$.
From the choice of $G$ and $G/K/N_i/K\simeq G/N_i\in s\mathfrak{U}$ for $i=1,2$
it follows $G/K/(N_1/K\cap N_2/K)\simeq G/N_1 \cap N_2\in s\mathfrak{U}$.
This contradicts the choice of $G$.

Let $N_1\cap N_2=1$. For every Sylow $p$-subgroup $P$ of $G$ a quotient group
$PN_i/N_i\in \text{\rm Syl}_p(G/N_i), i=1, 2$. From the strongly supersolubility of
$G/N_i$
it follows $PN_i/N_i$ is submodular in $G/N_i, i=1, 2$.
By 3) of Lemma 1.1,
$PN_i$ is submodular in $G, i=1, 2$. From properties of Sylow subgoups and 5)
of Lemma 1.1, it follows $PN_1\cap PN_2=P(N_1\cap N_2)=P$ is submodular in $G$.
So $G\in s\mathfrak{U}$. This contradicts the choice of $G$.
Statement 3) is proved.

Statement 4) follows from 3).

Prove Statement 5). Let $G/\Phi(G)\in s\mathfrak{U}$. From $s\mathfrak{U}\subseteq
\mathfrak{U}$ and the saturation of the class of groups
$\mathfrak{U}$ it follows $G\in \mathfrak{U}$. So
$F(G/\Phi(G))=F(G)/\Phi(G)$. Then $G/F(G)\simeq G/\Phi(G)/F(G/\Phi(G))
\in \mathfrak{B}$, i.e. $G\in s\mathfrak{U}$.

Statement 6) follows from 1)--3) and 5). Theorem is proved.

\medskip

{\bf Theorem 2.6.} {\sl The class of all strongly supersolubility groups is
a local formation and has a local screen $f$ such that $f(p)=\mathfrak{A}(p-1)\cap\mathfrak{B}$
for any prime $p$. }

{\bf Proof.} Since $f(p)=\mathfrak{A}(p-1)\cap\mathfrak{B}$ is a formation, $f$ is a local screen. Let a local formation $LF(f)$ be defined by a screen $f$. Let's show that $s\mathfrak{U} = LF(f)$.

Let $G\in s\mathfrak{U}$ and $H/K$ be any chief factor of $G$. From $G/F(G)\in
\mathfrak{B}$ and $F(G)\leq C_G(H/K)$ it follows $G/C_G(H/K)\in \mathfrak{B}$. Since $G$ is supersoluble, $|H/K|=p$ is some prime. By
Lemma 1.5,
$G/C_G(H/K)\in \mathfrak{A}(p-1)$. So $G/C_G(H/K)\in f(p)$. Then
$G\in LF(f)$ and $s\mathfrak{U}\subseteq LF(f)$.

Let now $G\in LF(f)$. Then $G/C_G(H/K) \in
f(p)\subseteq \mathfrak{B}$ for any chief factor $H/K$ of $G$ and $p\in \pi(H/K)$. Since $\mathfrak{B}$ is a formation, we conclude
$G/F(G) \in \mathfrak{B}$.
So $LF(f) \subseteq
s\mathfrak{U}$. Thus $LF(f) = s\mathfrak{U}$. Theorem is proved.

\medskip
{\bf Theorem 2.7.} {\sl Let the group $G=AB$ be the product of nilpotent 
subgroups $A$ and $B$. If $A$ and $B$ are submodular in $G$, then $G$ is strongly supersoluble.}

{\bf Proof.} Let $G$ be a counterexample of minimal order to Theorem.
Then, by Theorem of Wielandt-Kegel \cite{wiel,keg}, $G$ is soluble.
Let $N$ be a minimal normal subgroup of $G$.
Then $AN/N \simeq A/A \cap N \in \mathfrak{N},\ BN/N \simeq B/B \cap N \in \mathfrak{N}$.
  By 6) and 2) of Lemma 1.1, $AN/N$ and $BN/N$ are submodular in  $G/N$. By the choice of $G$, it follows $G/N \in s\mathfrak{U}$. By 3)--5) of Theorem 2.5, we conclude that $N$ is the only minimal normal subgroup in $G$ and $\Phi(G) = 1$. Then $G=MN$,
where $M$ is a maximal subgroup in $G$, $M \cap N = 1$, $N=C_G(N)$ and $|N|=p^n$
for some prime $p$. By 1) of Lemma 1.3, $A \cap B = 1$. From $N \subseteq A \cup B$
it follows either $N \leq A$, or $N \leq B$. Without loss of generality, we may suppose $N \leq A$. Then, by 3) of Lemma 1.3, $A$ is a $p$-subgroup, $B$ is a $p'$-subgroup.

Let $q$ be the largest prime divisor of $|G|$. If $q \neq p$,
 then $B$ has some Sylow $q$-subgroup $Q$ of $G$. From $Q \unlhd B$
and the submodularity of $B$ in $G$ it follows that $Q$ is submodular in $G$.
By Lemma 2.1, $Q \unlhd G$. Then $N \leq Q$. We get a contradiction with $q \neq p$.
So $q=p$. In view of Lemma 2.1, $A \unlhd G$. By Lemma 1.4, $O_p(M)=1$.
 Then $M \cap A = 1$ and $A=N$, $B$ is a maximal subgroup of $G$ and $B_G=1$.
Hence, $B$ is a maximal modular subgroup in $G$.
 In view of $B_G=1$ by Lemma 1.2, we conclude $|G|=pr$, where $r$ is a prime and $p\not= r$.
So $G \in s\mathfrak{U}$. This contradicts the choice of $G$. Theorem is proved.

\medskip

In Theorem 2.7 we can't discard the submodularity of one of subgroups.

{\bf Example 2.8.} In group $G=AB$, where $A\simeq Z_{17}$ and $B\simeq Aut(Z_{17}) \simeq Z_{16}$,
the subgroup $A$ is submodular, but the subgroup $B$ is not submodular in $G$. The group $G$ is supersoluble, but not strongly supersoluble. The example also shows that $s\mathfrak{U}\not=\mathfrak{U}$.

\medskip

{\bf Theorem 2.9.} {\sl A group $G$ is strongly supersoluble if and only if $G$ is metanilpotent and any Sylow subgroup of $G$ is submodular in $G$.}

{\bf Proof.} Necessity follows from that the strongly supersoluble group is supersoluble, and so it has a nilpotent commutator subgroup,
i.e. it is metanilpotent.

Sufficiency. Let there exists metanilpotent groups which have all Sylow subgroups are submodular in a group, but the group is not strongly supersoluble.
Let's choose from them a group $G$ of the smallest order. 
Let $N$ be a minimal normal subgroup of $G$. Then $G/N \in s\mathfrak{U}$
in view of the choice of $G$. Since, by Theorem 2.5, the class $s\mathfrak{U}$ is a saturated formation, then $N$ is the only minimal normal subgroup in $G$ and $\Phi(G)=1$. So $N = C_G(N)$ and $G=NM$, where $M$ is a maximal subgroups of $G$, $M \cap N =1$. From the metanilpotency of $G$
 and $N=F(G)$ it follows $G/N \simeq M \in \mathfrak{N}$. Let $p$ be the largest prime divisor of $|G|$.
Since $G$ is Ore dispersive, it follows that $N$ is contained in some Sylow $p$-subgroup of $G$. In view of $O_p(M)=1$, we conclude that $N \in \text{Syl}_p(G)$ and $M$ is a $p'$-group.
Let $S \in \text{Syl}_q(M)$.

If $G=SN$, then, by Theorem 2.7, $G$ is strongly supersoluble. This contradicts the choice of $G$.

Let $G\not=SN$ for every $S \in \text{Syl}_q(M)$. Denote $L=SN$.
Then $L$ is strongly supersoluble by the choice of $G$.
From $C_G(N)=N$ it follows $O_{p'}(L)=1$. Then $N=F_p(L)$. By Lemma 1.7 
and Theorem 2.6, we get $S\simeq L/F_p(L)\in
\mathfrak{A}(p-1)\cap\mathfrak{B}$. Hence and from the nilpotency of $M$
it follows that $M \in \mathfrak{A}(p-1)$. Since $N=F_p(G)$,
we conclude that $M\simeq G/F_p(G)\in \mathfrak{A}(p-1)$. By Lemmas 1.7 and 1.8, 
$G$ is supersoluble. By Definition 2.2, $G$ is strongly supersoluble.
This contradicts the choice of $G$. Theorem is proved.

\bigskip

{\parindent=0mm

\textbf{\bf {\large 3. Groups with submodular Sylow subgroups}}

}
\bigskip

Denote $sm\mathfrak{U}=(\ G\ |$ every Sylow subgroup of the group $G$ is
 submodular in $G$ ).

\medskip

{\bf Theorem 3.1.} {\sl Let $G$ be a group. Then the following hold:

$1)$ if $G \in sm\mathfrak{U}$ and $H \leq G$, then $H \in sm\mathfrak{U}$;

$2)$ if $G \in sm\mathfrak{U}$ and $N \unlhd G$, then $G/N \in
sm\mathfrak{U}$;

$3)$ if $N_i \unlhd G$ and $G/N_i \in sm\mathfrak{U}$, $i=1,2$, then $G/N_1 \cap N_2 \in
sm\mathfrak{U}$;

$4)$ if $H_i \in sm\mathfrak{U}$, 
$H_i \unlhd G$, $i=1,2$ and $H_1 \cap H_2 =1$,   
then $H_1 \times H_2 \in sm\mathfrak{U}$;

$5)$ if $G/\Phi(G) \in sm\mathfrak{U}$, then $G \in sm\mathfrak{U}$;

$6)$ the class of groups $sm\mathfrak{U}$ is a hereditary saturated formation.
}

{\bf Proof.} The validity of Statements 1) and 2) of Theorem follows from Statements 1), 2) and
6) of Lemma 1.1, in view of $H \cap G_p \in \text{Syl}_p(H)$ for some
$G_p \in \text{Syl}_p(G)$ and $R/N = G_qN/N \in \text{Syl}_q(G/N)$ for some $G_q \in \text{Syl}_q(G)$.

Prove Statement 3). Let $G$ be a group of the smallest order such that $G/N_i \in sm\mathfrak{U},\ N_i \unlhd G,\ i=1,2$, but
$G/N_1 \cap N_2 \not \in sm\mathfrak{U}$.


We can suppose that $N_1 \cap N_2 = 1$. Let $P \in \text{Syl}_p(G)$.
 Then $PN_i/N_i \in \text{Syl}_p(G/N_i),\ i=1,2$. So $PN_i/N$ is submodular in $G/N_i$.
 By 3) of Lemma 1.1, $PN_i$ is submodular in $G$. By the property of Sylow subgroups and Statement 5) of Lemma 1.1,
$PN_1 \cap PN_2 = P(N_1 \cap N_2)=P$ is submodular in $G$, i.e. $G \in sm\mathfrak{U}$.
This contradiction completes the proof of 3).

Statement 4) follows from 3).

Prove Statement 5). Let $G$ be a group of the smallest order such that
$G/\Phi(G) \in sm\mathfrak{U}$, but $G \not \in sm\mathfrak{U}$. Then $G$ is soluble in view of Corollary 2.1.1 and the solubility of $\Phi(G)$.
Let $N$ be a minimal normal subgroups of $G$.
From $\Phi(G)N/N \subseteq \Phi(G/N)$ and by Statement 2) of Theorem, we conclude that
$G/N/\Phi(G/N)
\in sm\mathfrak{U}$. Since $|G/N| < |G|$, $G/N \in sm\mathfrak{U}$. From Statement 3) it follows that $N$ is the only minimal subgroup of $G$, $|N|=p^n$ for some prime $p$
and $O_{p'}(G)=1$. Hence, $N \subseteq \Phi(G)$.

Let $Q \in \text{Syl}_q(G)$. From $QN/N \in \text{Syl}_q(G/N)$ it follows that $QN/N$ is submodular in $G/N$. By Statement 2) of Lemma 1.1, $QN$ is submodular in $G$.

If $p=q$, then $QN=Q$ is submodular in $G$.

Let $p \neq q$. Let's consider 2 cases:

(1) $|\pi(G)|=2$. Then $G/N = QN/N \cdot P/N$, where $P \in \text{Syl}_p(G)$.
By Theorem 2.7, $G/N$ is strongly supersoluble. Since the class $s\mathfrak{U}$ of all strongly supersoluble groups is a saturated formation by Theorem 2.5, then from
$G/\Phi(G) \simeq G/N/\Phi(G)/N \in s\mathfrak{U}$ it follows that $G \in s\mathfrak{U}
\subseteq sm\mathfrak{U}$. This contradicts the choice of $G$.
\medskip

(2) $|\pi(G)|>2$. Then $G/N \neq H/N = QN/N \cdot R/N$, where $R/N \in \text{Syl}_p(G/N)$.
By Statement 1) of Theorem $H/N \in sm\mathfrak{U}$. So Sylow $q$-subgroup $QN/N$ and Sylow $p$-subgroup $R/N$ of $H/N$ are submodular in $H/N$. By Theorem 2.7, $H/N$ is strongly supersoluble. Since by Theorem 2.6 $s\mathfrak{U}$ is a local formation and $H=QR$ is a Hall $\{p,q\}$-subgroup of $G$, then apply Corollary 16.2.3 from \cite{shem}.
 From $H\Phi(G)/\Phi(G) \simeq H/H \cap \Phi(G) \simeq H/N/H \cap \Phi(G)/N \in s\mathfrak{U}$
it follows $H \in s\mathfrak{U}$. By Statement 1) of Theorem 2.5, $QN \in s\mathfrak{U}$.
Then $Q$ is submodular in $QN$, and so, $Q$ is submodular in $G$. The arbitrariness of the choice of
$Q$ implies $G \in sm\mathfrak{U}$. This contradiction completes the proof of Statement 5).

Statement 6) follows from Statements 1)-3) and 5). Theorem is proved.

\medskip


{\bf Lemma 3.2.} {\sl The following statements are held.

$1)$ The class of groups $(G \in \mathfrak{S}\ |\ \text{\rm Syl}(G) \subseteq
\mathfrak{B})$ is a hereditary formation.

$2)$ For any prime $p$ the class of groups
$(G \in \mathfrak{S}\ |\ \text{\rm Syl}(G) \subseteq \mathfrak{A}(p-1)\cap\mathfrak{B})$
is a hereditary formation.}

{\bf Proof.} Prove Statement 1). Denote $\mathfrak{H}=(G \in \mathfrak{S}\ |\ \text{\rm Syl}(G) \subseteq
\mathfrak{B})$. Clearly, if
$H \leq G \in \mathfrak{H}$ and $N \unlhd G$, then $H \in \mathfrak{H}$ and $G/N \in \mathfrak{H}$.

Let's show by induction on $|G|$, if
$N_i \unlhd G$ and  
$G/N_i \in \mathfrak{H},\ i=1,2$,
then $G/N_1 \cap N_2 \in \mathfrak{H}$. If $K = N_1 \cap N_2 \neq 1$,
then from $|G/K| < G$ and $G/K/N_i/K \simeq G/N_i \in \mathfrak{H}$ it follows  $G/K/N_1/K \cap N_2/K \simeq G/N_1 \cap N_2 \in \mathfrak{H}$.
Let $N_1 \cap N_2 = 1$. Let $P \in \text{Syl}_p(G)$.
From $G/N_i \in \mathfrak{H}$ it follows that $PN_i/N_i \simeq P/P \cap N_i$
is an elementary abelian group. Since the class of all abelian groups  $\mathfrak{A}$ is a formation, then $P/(P \cap N_1) \cap (P \cap N_2) \simeq P \in \mathfrak{A}$.
We will show that $P$ is an elementary abelian $p$-group. Let $z \in P$, $|z|=p^n$
and $Z = \langle z \rangle$. From $ZN_i/N_i \leq  PN_i/N_i$ it follows that
$|ZN_i/N_i|=|Z/Z \cap N_i| \leqslant p$. Since $N_1 \cap N_2 = 1$,
then there exists $i \in \{1,2\}$ such that $Z \cap N_i = 1$. Then $|Z| = p^n,\ n=1$.
Since $P$ is a direct product of cyclic subgroups, we get $P\in \mathfrak{B}$. So $G \in \mathfrak{H}$.

Statement 2) is being proved similarly taking into account
$\mathfrak{A}(p-1)\cap\mathfrak{B}$ is a hereditary formation.
Lemma is proved.

\medskip

{\bf Lemma 3.3.} {\sl A local formation $LF(f)$ with a local screen $f$ such that
$f(p)=(H \in \mathfrak{S}\ |\ \text{\rm Syl}(H) \subseteq
\mathfrak{A}(p-1)\cap\mathfrak{B})$
for any prime $p$, 
is a hereditary saturated formation.}

{\bf Proof.} By Theorem 1.6, $LF(f)$ is a saturated formation.

Let's prove the heredity of $LF(f)$. Let $G \in LF(f)$ и $R \leq G$.
Then $G$ has a chief series $1=G_0 < G_1 < \ldots < G_{n-1} < G_n =
G$ such that $G/C_G(G_i/G_{i-1}) \in f(p)$ for every $p \in \pi(G/G_{i-1})$ and
$i=1,\ldots, n$.
Let $R_{i-1}=R \cap G_{i-1},\ i=1,\ldots,n+1.$ Let $C_i=C_G(G_i/G_{i-1})$
 and $C_i^*=C_R(R_i/R_{i-1}),\ i=1,\ldots,n$. It is easy to see that
$R \cap C_i \leq C_i^*$.
From $RC_i/C_i \leq G/C_i \in f(p)$ and the heredity of $f(p)$ it follows that $RC_i/C_i \simeq R/R \cap C_i \in f(p)$. Then
$R/C_i^* \simeq R/R \cap C_i / C_i^*/R \cap C_i \in f(p)$.
Hence $R/C_R(H/K) \in f(p)$ for every chief factor $H/K$ of $R$ and $p \in \pi(H/K)$. So $R\in LF(f)$.
Lemma is proved.

\medskip

{\bf Theorem 3.4.} {\sl Every minimal non $sm\mathfrak{U}$-group is biprimary minimal non $s\mathfrak{U}$-group.}

{\bf Proof.} Let $G \in {\cal{M}}(sm\mathfrak{U})$ and $q$ be the smallest prime divisor of $|G|$. Every subgroup $H$ of a group $G$ belongs
$sm\mathfrak{U}$. By Corollary 2.1.1, $H$ is Ore disperive. So $H$ is
$q$-nilpotent. Let's consider two cases.

1) $G$ is $q$-nilpotent. Then $G=QP$, where $Q \in \text{Syl}_q(G),\ P \unlhd G$ and
$P$ is a Hall $q'$-subgroup of $G$. From $P \in sm\mathfrak{U}$ it follows the solubility of $P$. Then from $G/P \simeq Q$ we get the solubility of $G$.

Suppose that $\Phi(G)=1$. Let $N$ be a minimal normal subgroup of $G$.
Then $|N|=p^n$ for some $p \in \pi(G)$. Let $G=NM$, where $M$ is a maximal subgroup of $G$. In view of $G \in {\cal{M}}(sm\mathfrak{U})$ and
$G/N \simeq M / M \cap N \in sm\mathfrak{U}$, $N$ is the only minimal normal subgroup of $G$. Let $R$ is an arbitrary Sylow $r$-subgroup of $G$. From $G/N \in sm\mathfrak{U}$ we conclude that $RN/N$ is submodular in $G/N$. By 3) of Lemma 1.1, $RN$ is submodular in $G$. If $RN \neq G$,
then from $RN \in sm\mathfrak{U}$ we get $R$ is submodular in $G$.
This contradicts with $G \not \in sm\mathfrak{U}$. Hence $RN = G$ is a biprimary group. Since every subgroup $T$ of $G$ belongs
$sm\mathfrak{U}$, then, by Theorem 2.7, $T \in s\mathfrak{U}$. From
$s\mathfrak{U} \subseteq sm\mathfrak{U}$ it follows that $G \not \in s\mathfrak{U}$,
i.e. $G \in {\cal{M}}(s\mathfrak{U})$.

Let $\Phi(G) \neq 1$. Then $\Phi(G/\Phi(G)) = 1$. Since $sm\mathfrak{U}$ is saturated, it follows $G/\Phi(G) \not \in sm\mathfrak{U}$. Then
$G/\Phi(G) \in {\cal{M}}(sm\mathfrak{U})$. As proved above, $G/\Phi(G)$ is a biprimary group and $G/\Phi(G) \not \in s\mathfrak{U}$. Hence $G \in {\cal{M}}(s\mathfrak{U})$.

2) $G$ is not $q$-nilpotent. By Theorem 5.4 of \cite[гл. IV]{Hup}, $G$ is a Schmidt group. Since every subgroup $T$ of $G$ is nilpotent, then $T \in s\mathfrak{U}$.
Then $G \in {\cal{M}}(s\mathfrak{U})$. Theorem is proved.

\medskip

Recall that a subgroup $H$ of a group $G$ is called
{\it ${\rm K}$-$\Bbb{P}$-subnormal \emph{\cite{vasTyt}} in $G$} ,
if there exists a chain of subgroups
\begin{equation}
H = H_0\leq H_1\leq\cdots\leq H_{n-1}\leq H_n=G
\label{1}
\end{equation}
such that either $H_{i-1}$ is normal in $H_{i}$, or $|H_{i} : H_{i-1}|$ is a prime
for every $i = 1,\ldots, n$. 

 If $H=G$ or in a chain~\eqref{1} the index $|H_{i} : H_{i-1}|$ is a prime for every
$i = 1,\ldots, n$, then $H$ is called {\it $\Bbb{P}$-subnormal in $G$} \cite{vas}.

\medskip

{\bf Lemma 3.5.} {\sl Let $H$ be a submodular Sylow subgroup of a group $G$.
Then the following conditions are held:

$1)$ $H$ is ${\rm K}$-$\Bbb{P}$-subnormal in $G$;

$2)$ if $G$ is soluble, then $H$ is $\mathbb{P}$-subnormal in $G$.}

\medskip

{\bf Proof.} Prove 1) by induction on $|G|$. We can suppose that $H\not=G$. Then $H$ is contained in a maximal modular subgroup $M$ of $G$. By 1) of Lemma 1.1
 and $|M| < |G|$, it follows that $H$ is ${\rm K}$-$\Bbb{P}$-subnormal in $M$.
By Lemma 1.2 either $M$ is normal in $G$, or $G/M_G$ is non-abelian of order $pq$,
where $p$ and $q$ are primes. Hence, if $M_G\not= M$ we have $|G : M|=
|G/M_G: M/M_G|$ is a prime. This means that $M$ is ${\rm
K}$-$\Bbb{P}$-subnormal in $G$. So $H$ is ${\rm
K}$-$\Bbb{P}$-subnormal in $G$.

Statement 2) follows from 1), since in a soluble group
${\rm K}$-$\Bbb{P}$-subnormal subgroup is
$\Bbb{P}$-subnormal. Lemma is proved.

\medskip

By Lemma 3.5, it follows that $sm\mathfrak{U}\subseteq \text{w}\mathfrak{U}$,
where $\text{w}\mathfrak{U}$ is the class of all groups with
$\mathbb{P}$-subnormal Sylow subgroups. Example 2.8 shows that $sm\mathfrak{U}\not= \text{w}\mathfrak{U}$.

\medskip

{\bf Theorem 3.6.} {\sl The class of all groups with submodular Sylow subgroups is a local formation and has a local screen $f$ such that
$f(p)=(G \in \mathfrak{S}\ |\  \text{\rm Syl}(G) \subseteq \mathfrak{A}(p-1)\cap\mathfrak{B})$
for every prime $p$.}

{\bf Proof.} Since $f(p)=(G \in \mathfrak{S}\
|\  \text{\rm Syl}(G) \subseteq \mathfrak{A}(p-1)\cap\mathfrak{B})$
is a formation, $f$ is a local screen. Let a local formation $LF(f)$ be defined by a screen $f$.
Denote $\mathfrak{F}=LF(f)$. By Theorem 2.10 \cite{vas}, the class of groups
$\text{w}\mathfrak{U}$
is a local formation and has a local screen $h$ such that
$h(p) = (G \in \mathfrak{S}\
|\  \text{\rm Syl}(G) \subseteq \mathfrak{A}(p-1))$ for every prime $p$. Hence $\mathfrak{F}\subseteq\text{w}\mathfrak{U}$. In view of Proposition 2.8
\cite{vas}, $\mathfrak{F}$ consists of Ore dispersive groups.

Show that $\mathfrak{F} \subseteq sm\mathfrak{U}$. Let $G$ be a group of the smallest order from $\mathfrak{F} \backslash sm\mathfrak{U}$. By Lemma 3.3, $\mathfrak{F}$ is a hereditary formation. Hence $G$ is a minimal non $sm\mathfrak{U}$-group.
Since $\mathfrak{F} \subseteq \mathfrak{S}$ and $sm\mathfrak{U}$ is a hereditary formation by 6) of Theorem 3.1, then $G$ has the unique minimal normal subgroup $N$, $N=C_G(N)$ is an elementary abelian $p$-subgroup for some prime $p$, $\Phi(G)=1$. Then $G = NM$, where $M$ is a maximal subgroup of $G$.
By Lemma 1.4, $O_p(M) = 1$. Since $G$ is Ore dispersive, it follows that
$N \in \text{Syl}_p(G)$
 and $p$ is the largest prime divisor of $|G|$. By Theorem 3.4, $G$ is a biprimary minimal non $s\mathfrak{U}$-group. Hence and from $G/C_G(N) \simeq M \in f(p)$
we get that $M$ is an elementary abelian $q$-group and $q$ divides $(p-1)$. Since
 $M$ is $\mathbb{P}$-subnormal in $G$ it follows that $|G:M|=p$ and $|N|=p$.
Then $M \simeq G/N$ is isomorphically embedded in a cyclic group of order $p-1$.
So $|M|=q$. Hence $G \in s\mathfrak{U} \subseteq sm\mathfrak{U}$.
This contradicts the choice of $G$. So $\mathfrak{F} \subseteq sm\mathfrak{U}$.

Prove that $sm\mathfrak{U} \subseteq \mathfrak{F}$. Let $G$ be a group of the smallest order from $sm\mathfrak{U} \backslash \mathfrak{F}$. Since $G \in sm\mathfrak{U}$, then
$G$ is soluble by Corollary 2.1.1. Since $sm\mathfrak{U}$ and $\mathfrak{F}$ are saturated formations, then $\Phi(G)=1$. In $G$ there exists the unique minimal normal subgroup $N = C_G(N) = F(G)$, $|N|=p^n$ for some $p \in \pi(G)$.
Then $G=NM$, where $M$ is a maximal subgroup in $G$, $N \cap M = 1$.
By Corollary 2.1.1, $G$ is Ore dispersive.
 Then $p$ is the largest prime divisor of $|G|$. In view of Lemma 1.4, $N \in \text{Syl}_p(G)$
 and $M$ is a $p'$-groups. Let $Q \in \text{Syl}_q(M)$. Then $Q \in \text{Syl}_q(G)$.

Suppose that $QN \neq G$. In view of $QN \in sm\mathfrak{U}$ and by Theorem 2.7,
$QN$ is strongly supersoluble.
From $N \leq F(QN)$, $N=C_G(N)$ and $O_{p'}(QN)=1$, it follows that $N = F_p(QN)$.
By Lemma 1.7, 
$Q \simeq QN/F_p(QN) \in
\mathfrak{A}(p-1)\cap\mathfrak{B}$.
Then $M \in f(p)$ in view of the arbitrariness of the choice of $Q$.
By Lemma 1.7 $G \in \mathfrak{F}$.  
This contradicts the choice of $G$. 

Hence $QN=G$. By Theorem 2.7, $G$ is strongly supersoluble. By Lemma 2.6,
$G/C_G(N)=G/N\in \mathfrak{A}(p-1)\cap\mathfrak{B}$. So
$G \in \mathfrak{F}$. This contradicts the choice of $G$. Theorem is proved.

\medskip


{\bf Theorem 3.7.} {\sl Let $G$ be a group in which every Sylow subgroup is submodular in $G$. Then the following conditions are held:

$1)$ every metanilpotent subgroup of $G$ is strongly supersoluble;

$2)$ every biprimary subgroup of $G$ is strongly supersoluble;

$3)$ if $N$ is the smallest normal subgroups of $G$ such that $G/N$ is a group with elementary abelian Sylow subgroups, then $N$ is nilpotent.

}

{\bf Proof.} Statement 1) follows from the heredity of the class of groups
$sm\mathfrak{U}$ and Theorem~3.6.

Statement 2) follows from the heredity of the class of groups $sm\mathfrak{U}$ and Theorem 2.9.

Prove Statement 3). By Lemma 3.2, the class of groups
$\mathfrak{H}=(G \ |\ \text{\rm Syl}(G) \subseteq
\mathfrak{B})$ is a hereditary formation.
 By Statement 10 of \cite[с.~36]{shem}, $\mathfrak{N}\mathfrak{H}$ has a local screen $f$ such that
$f(p)=\mathfrak{H}$ for every prime $p$. In view 
of Theorem 3.6, 
$sm\mathfrak{U} \subseteq \mathfrak{N}\mathfrak{H}$. Then $N=
G^{\mathfrak{H}}\in \mathfrak{N}$. Theorem is proved.

{\bf Theorem 3.8.} {\sl Every Sylow subgroup of the group $G$ is submodular 
in $G$ if and only if the group $G$ is Ore dispersive and every its biprimary 
subgroup is strongly supersoluble.}

{\bf Proof.} Let the class of groups $\mathfrak{F}=(G\ |$ the group $G$
is Ore dispersive and every biprimary subgroup of $G$ is strongly supersoluble).

 If $G \in sm\mathfrak{U}$, then $G \in \mathfrak{F}$ in view of Corollary 2.1.1 and by 2) of Theorem 3.7.

Let $G$ be a group of the smallest order belonging $\mathfrak{F}
\backslash sm\mathfrak{U}$. Since $s\mathfrak{U} \subseteq sm\mathfrak{U}$,
then $|\pi(G)|>2$. Since $G$ is Ore dispersive, then for the largest prime
$p \in \pi(G)$ and $P \in \text{Syl}_p(G)$, $P$ is normal in $G$.
Since $G$ is soluble, there exists a Hall $p'$-subgroup $H$ of $G$.
From $H \neq G$ and $H \in \mathfrak{F}$ if follows that $H \in sm\mathfrak{U}$.
Note that $|\pi(H)| \geqslant 2$. Let $Q \in \text{Syl}_q(G),\ q \neq p$.
Since $QP/P \in \text{Syl}_q(G/P)$ and $G/P \simeq H \in sm\mathfrak{U}$,
then $QP/P$ is submodular in $G/P$. By 3) of Lemma 1.1, $QP$ is submodular in $G$. From $QP \in s\mathfrak{U}$, it follows that $Q$
is submodular in $QP$. So $Q$ is submodular in $G$. This means that
$G \in sm\mathfrak{U}$. We got a contradiction with the choice of $G$. Theorem is proved.

\renewcommand{\refname}{\hfill\normalsize\bf References\hfill}

V.\,A.\,Vasilyev

Francisk Skorina Gomel State University,
Sovetskaya str., 104,
Gomel 246019,
Belarus.
E-mail address: vovichx@mail.ru


\end{document}